\documentclass[11pt]{article}

\def\init{\setcounter{equation}{0}}

\newtheorem{theorem}{Theorem}[section]

\newtheorem{lemma}[theorem]{Lemma}

\newtheorem{corollary}[theorem]{Corollary}

\newcommand{\qed}{{\vrule height 1ex width 1ex depth -.1ex }}

\def\nn{{\bf N}}
\def\zz{{\bf Z}}

\def\cala{{\cal A}}
\def\calb{{\cal B}}

\def\proof{\noindent{\bf Proof. }}

\def\mod{\hbox{mod}\,}

\newcommand{\beq}{\begin{equation}}
\newcommand{\eeq}{\end{equation}}

\def\aa{\alpha}

\begin{document}

\title{A characterization of finite sets that tile the integers}
\author{Andrew Granville, Izabella {\L}aba, and Yang Wang}
\maketitle


\section{Introduction}
\label{intro}
\init


A set of integers $A$ is said to {\em tile the integers} 
if there is a set $C\subset\zz$ such that every integer $n$ can be written
in a unique way as $n=a+c$ with $a\in A$ and $c\in C$. Throughout this
paper we will assume that $A$ is finite.  It is well known (see
\cite{New}) that any tiling of $\zz$ by a finite set $A$ must be periodic:
$C=B+M\zz$ for some finite set $B\subset \zz$ such that $|A|\,|B|=M$.  We
then write $A\oplus B=\zz/M\zz$.  

Newman \cite{New} gave a characterization of all sets $A$ which tile
the integers and such that $|A|$ is a prime power. Coven and Meyerowitz
\cite{CM} found necessary and sufficient conditions for $A$ to tile $\zz$
if $|A|$ has at most two prime factors.  To state their result we need to
introduce some notation. Without loss of generality we may assume that
$A,B\subset\{0,1,\dots\}$ and that $0\in A\cap B$. Define the characteristic 
polynomials
\[
A(x)=\sum_{a\in A}x^a,\ B(x)=\sum_{b\in B}x^b .
\]
Then $A\oplus B=\zz / M\zz$ is equivalent to
\begin{equation}\label{d.e1}
 A(x)B(x)=1+x+\dots+x^{M-1}\ (\mod (x^M-1)).
\end{equation}
Let $\Phi_s(x)$ denote the $s$-th cyclotomic
polynomial, i.e. the monic, irreducible polynomial whose roots are the
primitive $s$-th roots of unity. We then have $x^n-1=\prod_{s|n}\Phi_s(x)$,
and (\ref{d.e1}) holds if and only if
\begin{equation}\label{d.e2}
  |A||B|=M\hbox{ and }\Phi_s(x)\ |\ A(x)B(x)\hbox{ for all }s|M,\ s\neq 1.
\end{equation}
Let $S_A$ be the set of prime powers
$p^\alpha$ such that $\Phi_{p^\alpha}(x)$ divides $A(x)$.  Then the
Coven-Meyerowitz conditions are:

\medskip
{\it (T1) $A(1)=\prod_{s\in S_A}\Phi_s(1)$,}

\medskip
{\it (T2) if $s_1,\dots,s_k\in S_A$ are powers of different
primes, then $\Phi_{s_1\dots s_k}(x)$ divides $A(x)$.}
\medskip

\noindent It is proved in \cite{CM} that:
\begin{itemize}
\item if $A$ satisfies (T1), (T2), then it tiles $\zz$;
\item if $A$ tiles $\zz$ then (T1) holds;
\item if $A$ tiles $\zz$ and $|A|$ has at most two prime factors,
then (T2) holds.
\end{itemize}
The first two statements are relatively simple to prove and hold regardless
of the size of $A$; the main difficulty is in proving the third one. 
The proof given by Coven and Meyerowitz relies crucially on a result of Sands
\cite{S1}: if $A\oplus B=\zz/M\zz$ and $M$ has at most two prime divisors,
then one of $A$, $B$ must be contained in $m\zz$ for some $m|M$, $m\neq 1$.  A theorem 
of Tijdeman \cite{Tij} implies that if $A$ tiles the integers, then there
exists a tiling $A\oplus B$ such that $|B|$ has the same prime factors as $|A|$.
Therefore if $|A|$ has at most two prime factors, there is a tiling to which
Sands' result applies. The authors then decompose this tiling and proceed by induction
in $|A|$.  

It seems very hard to verify whether (T2) holds for all sets which tile the 
integers.  There is no analogue of Sands' result if $M$ has three or more
prime factors, as shown in \cite{Sz}, \cite{LS}; hence the methods of Coven 
and Meyerowitz do not extend to more general sets.
The purpose of this paper is to settle,
for the first time, a three-prime case. 

\begin{theorem}
Let $A,B$ be two sets of integers such that $|A|=p^\alpha q^\beta r^\gamma$ and
$|B|=pqr$, where $p,q,r$ are distinct primes. Assume that $A\oplus B=\zz/M\zz$,
where $M=|A|\,|B|$.  If $\Phi_p(x),\Phi_q(x),\Phi_r(x)$ divide $A(x)$, then so 
do $\Phi_{pq}(x), \Phi_{pr}(x), \Phi_{qr}(x), \Phi_{pqr}(x)$. \label{three}
\end{theorem}

Equivalently, if the elements of $A$ are equi-distributed modulo
$p$, $q$, and $r$, then they are also equi-distributed modulo $pqr$. Observe
that this reformulation of (T2) does not require
the elements of $A$ to be nonnegative.

We remark that by the results of \cite{LW1}, \cite{LW2}, \cite{L},
proving (T2) for all finite sets which tile the integers would essentially 
resolve one part of Fuglede's spectral set conjecture \cite{Fug}
in dimension 1.

Our main tool in proving Theorem \ref{three} is the following identity.  

\begin{theorem}\label{exponentials}
For any finite $A,B\subset \zz$, let
\[
A_m=\#\{(a,a')\in A\times A:\ (a-a',N)=m\},\
B_m=\#\{(b,b')\in B\times B:\ (b-b',N)=m\}.
\]
Then
\begin{equation}
\sum_{m|N} \frac{A_mB_m}{\phi(N/m)}
=\frac{1}{N}\sum_{d|N}\frac{\cala_d\calb_d}{\phi(d)},
\label{main.1}
\end{equation}
where
\[
\cala_d=\sum_{\xi:\Phi_d(\xi)=0}|A(\xi)|^2,\
\calb_d=\sum_{\xi:\Phi_d(\xi)=0}|B(\xi)|^2.
\]
\end{theorem}

Here, as usual, $\phi(n)$ is the 
Euler function and $(m,n)$ denotes the greatest common divisor of
$m$ and $n$.  We adopt the convention that $(n,0)=n$ for any $n\neq 0$.

We also observe that Theorem \ref{exponentials} extends the following result 
of Sands \cite{S1}. 

\begin{theorem}\label{sands}\cite{S1}
Let $A,B$ be two subsets of $\zz$ such that the elements of each of them
are distinct modulo $M$. Define $D_A=\{(a-a',M):\ a,a'\in A,\ a\neq a'\}$ and 
$D_B=\{(b-b',M):\ b,b'\in B,\ b\neq b'\}$. Then $A\oplus B=\zz/M\zz$ if 
and only if $|A|\,|B|=M$ and $D_A\cap D_B=\emptyset$.  
\end{theorem}

Our Theorem \ref{exponentials} provides an alternative proof of Theorem
\ref{sands}; furthermore, it implies Theorem \ref{divisors} below.

\begin{theorem}\label{divisors}
Define $A,B,D_A,D_B$ as in Theorem \ref{sands}.
If $D_A\cap D_B=\emptyset$, then $|A|\,|B|\leq M$; the equality
holds if and only if $A\oplus B=\zz/M\zz$.
\end{theorem}


\section{Proof of Theorems \ref{exponentials}--\ref{divisors}}
\label{identity}
\init


\noindent{\bf Proof of Theorem \ref{exponentials}.}  Fix $A,B\subset\zz$
and $N\in\nn$.  As usual, $\mu(n)$ is the M\"obius function and
$e(t)=e^{2\pi it}$. Let $d|N$, then for any $t\in\zz$
\begin{equation}\label{a.e1}
 \frac{1}{d}\sum_{\begin{array}{c}
 j=0\\ \frac{N}{d}\,|\,j\end{array}}^{N-1}
 e\Big(\frac{tj}{N}\Big)=\left\{
 \begin{array}{ll}1\ &\hbox{if }d\ |\ t,\\[3mm]0\ &\hbox{if }d\not{|}\ t.
 \end{array}\right.
\end{equation}
Let $\chi_I$ denote the characteristic function of the set $I$. Then for
$m|N$,
\[ \chi_{(N,t)=m}=\chi_{m|t}\,\chi_{(\frac{N}{m},\frac{t}{m})=1}
  =\chi_{m|t}\sum_{l|\frac{N}{m},\ l|\frac{t}{m}}\mu(l)\]
\[=\sum_{l|\frac{N}{m}}\mu(l)\chi_{lm|t}\]
\[=\sum_{l|\frac{N}{m}}\frac{\mu(l)}{lm}\sum_{\begin{array}{c}
  j=0\\ \frac{N}{lm}|j\end{array}}^{N-1} e\Big(\frac{tj}{N}\Big)
\]
\begin{equation}\label{a.e2}
 =\sum_{j=0}^{N-1}e\Big(\frac{tj}{N}\Big)  \sum_{l|\frac{N}{m},\ \frac{N}{lm}|j}
 \frac{\mu(l)}{lm}, 
\end{equation}
where we used (\ref{a.e1}) and that $\sum_{v|u}\mu(v)=\chi_{u=1}$.
Taking $v=N/m$, we deduce that
\[ A_m=\#\{(a,a')\in A\times A:\ (a-a',N)=m\}
=\sum_{a,a'\in A}\chi_{(a-a',N)=N/v}\]
\[=\sum_{a,a'\in A}\sum_{j=0}^{N-1}e\Big(\frac{(a-a')j}{N}\Big) 
\sum_{l|v,\ \frac{v}{l}|j}\frac{\mu(l)}{l}\]
\begin{equation}\label{a.e3}
  =\frac{1}{N}\sum_{j=0}^{N-1}\Big|\sum_{a\in A}e\Big(\frac{aj}{N}\Big) \Big|^2
  \sum_{d|v,\ d|j}d\mu(v/d),
\end{equation}
where we substituted $d=v/l$.  Let 
\[s_J=\Big|A\Big(e\Big(\frac{J}{N}\Big)\Big)\Big|=\Big|\sum_{a\in A}e\Big(\frac{aJ}{N}\Big)\Big|,\ 
 t_J=\Big|B\Big(e\Big(\frac{J}{N}\Big)\Big)\Big|=\Big|\sum_{b\in B}e\Big(\frac{bJ}{N}\Big)\Big|,\]
then from (\ref{a.e3}) we have
\[ \sum_{m|N} \frac{A_m\,B_m}{\phi(N/m)}=\sum_{v|N}\frac{1}{\phi(v)}
\left(\frac{1}{N}\sum_{J=0}^{N-1}s_J^2\sum_{d|v,\ d|J}d\mu(v/d)\right)
\left(\frac{1}{N}\sum_{I=0}^{N-1}t_I^2\sum_{e|v,\ e|I}e\mu(v/e)\right)\]
\begin{equation}\label{a.e4}
  =\frac{1}{N^2}\sum_{I,J=0}^{N-1} s_J^2\, t_I^2
  \left(\sum_{v|N}\frac{1}{\phi(v)}\sum_{d|v,\ d|J}d\mu(v/d)\sum_{e|v,\ e|I}e\mu(v/e)\right).
\end{equation}

Let $g=(I,J,N)$, $r=(J,N)/g$, and $s=(I,N)/g$ so that $(r,s)=1$. Then
\[ \sum_{v|N}\frac{1}{\phi(v)}\sum_{d|v,\ d|J}d\mu(v/d)\sum_{e|v,\ e|I}e\mu(v/e)
=\sum_{v|N}\frac{1}{\phi(v)}\sum_{d|(v,rg)}d\mu(v/d)\sum_{e|(v,sg)}e\mu(v/e)\]
\begin{equation}\label{a.e5}
  =\prod_{p^\alpha||N}\left(
  \sum_{i=0}^\alpha \frac{1}{\phi(p^i)}\sum_{d|(p^i,rg)}d\mu\Big(\frac{p^i}{d}\Big)
  \sum_{e|(p^i,sg)}e\mu\Big(\frac{p^i}{e}\Big)\right),
\end{equation}
since all the functions involved are multiplicative. Now 
\[ \sum_{d|(p^i,t)}d\mu\Big(\frac{p^i}{d}\Big)=\left\{
 \begin{array}{llll}
 1 &\hbox{ if }i=0,&&\\[3mm]
 p^i-p^{i-1} &\hbox{ if }i\geq 1&\hbox{ and }&p^i\ |\ t,\\[3mm]
 -p^{i-1} &\hbox{ if }i\geq 1&\hbox{ and }&p^{i-1}\ ||\ t,\\[3mm]
 0 &\hbox{ if }i\geq 1&\hbox{ and }&p^i\ \not|\ t.
 \end{array}\right. \]
Write $p^\gamma||g$ and $p^\delta||rs$ so that $\gamma+\delta\leq\alpha$.
Therefore
\[\sum_{i=0}^\alpha \frac{1}{\phi(p^i)}\sum_{d|(p^i,rg)}d\mu\Big(\frac{p^i}{d}\Big)
  \sum_{e|(p^i,sg)}e\mu\Big(\frac{p^i}{e}\Big)\]
\[ =1+\sum_{i=1}^\gamma \frac{1}{\phi(p^i)}(p^i-p^{i-1})^2+\left\{
 \begin{array}{lll}
 0&\hbox{ if }\gamma=a,&\\[3mm]
 \frac{1}{\phi(p^{\gamma+1})}(-p^\gamma)^2&\hbox{ if }\gamma<\alpha
 &\hbox{ and }\delta=0,\\[3mm]
 \frac{1}{\phi(p^{\gamma+1})}(-p^\gamma)(p^{\gamma+1}-p^\gamma)&\hbox{ if }\gamma<\alpha
 &\hbox{ and }\delta\geq 1\\[3mm]
 \end{array}\right. \]
\[ = \left\{\begin{array}{lll}
 p^\gamma &\hbox{ if }\gamma=a&\hbox{ (hence }\delta=0),\\[3mm]
 p^{\gamma+1}/(p-1)&\hbox{ if }\gamma<\alpha &\hbox{ and }\delta=0,\\[3mm]
 0&\hbox{ if }\gamma<\alpha &\hbox{ and }\delta\geq 1.\\[3mm]
 \end{array}\right. \]
We thus have a non-zero term in (\ref{a.e5}) if and only if $\delta=0$ for all $p$,
that is $r=s=1$, in other words $(I,N)=(J,N)=g$. In this case our answer is
$p^\gamma \,p^{\alpha-\gamma}/\phi(p^{\alpha-\gamma})=p^{\alpha}/\phi(p^{\alpha-\gamma})$.
Therefore (\ref{a.e5}) becomes
\[ \prod_{p^\alpha||N}\frac{p^\alpha}{\phi(p^{\alpha-\gamma})}=\frac{N}{\phi(N/g)}.\]
Substituting this into (\ref{a.e4}) gives
\begin{equation}\label{a.e6}
 \sum_{m|N}\frac{A_m\,B_m}{\phi(N/m)}=\frac{1}{N}\sum_{g|N}\frac{1}{\phi(N/g)}
\Big(\sum_{\begin{array}{c}\scriptstyle{I=0}\\ \scriptstyle{(I,N)=g}\end{array}}^{N-1}s_I^2\Big)
\Big(\sum_{\begin{array}{c}\scriptstyle{J=0}\\ \scriptstyle{(J,N)=g}\end{array}}^{N-1}t_J^2\Big).
\end{equation}
Let $N/g=d$ and $I=gi$, so that $(i,N/g)=(i,d)=1$. Then
\[\sum_{\begin{array}{c}\scriptstyle{I=0}\\ \scriptstyle{(I,N)=g}\end{array}}^{N-1}s_I^2
=\sum_{\begin{array}{c}\scriptstyle{i=0}\\ \scriptstyle{(i,d)=1}\end{array}}^{N-1}
\Big|A(e(ig/N))\Big|^2
=\sum_{\begin{array}{c}\scriptstyle{i=0}\\ \scriptstyle{(i,d)=1}\end{array}}^{N-1}
\Big|A(e(i/d))\Big|^2=\cala_d,\]
and similarly for $\sum t_J^2$.  Hence the right side of (\ref{a.e6}) equals
\[\frac{1}{N}\sum_{d|N}\frac{1}{\phi(d)}\cala_d\,\calb_d.\]
The theorem follows. \qed

\bigskip
{\bf Proof of Theorems \ref{sands} and \ref{divisors}.} Apply Theorem \ref{exponentials}
with $M=N$.  The term on the right side of (\ref{main.1}) with $d=1$ is $|A|^2\,|B|^2/M$,
and, since all elements of $A$ and $B$ are distinct modulo $M$, the term
on the left side of (\ref{main.1}) with $m=M$ is $|A|\,|B|=M$.  In particular, the left side
of (\ref{main.1}) is $\geq M$, since all the remaining terms are nonnegative.

We first deduce Theorem \ref{sands}.  We have $A\oplus B=\zz/M\zz$ if and
only if $|A|\,|B|=M$ and $\Phi_d(x)$ divides $A(x)$ or $B(x)$ for all $d|M$, $d\neq 1$.
This in turn is equivalent to $|A|\,|B|=M$ and $\cala_d\,\calb_d=0$ for
all $d|M$, $d\neq 1$.  Thus
$A\oplus B=\zz/M\zz$ if and only if $|A|\,|B|=M$ and the right side of 
(\ref{main.1}) equals $|A|^2\,|B|^2/M=M$.  But the left side of (\ref{main.1})
equals $M$ if and only if $A_m\,B_m=0$ for all $m|M$,
$m\neq M$, which in turn is equivalent to $D_A\cap D_B=\emptyset$.

Assume now that $D_A\cap D_B=\emptyset$. Then the left side of (\ref{main.1})
equals $M$, therefore so does the right side.  Using that the $d=1$ term
is $|A|^2\,|B|^2/M$ and that all other terms are nonnegative, we find that
$M\geq |A|^2\,|B|^2/M$, hence $|A|\,|B|\leq M$ and equality holds if and
only if all the terms with $d>1$ on the right are zero.  As above, the latter
together with the equality $|A|\,|B|=M$
is equivalent to $A\oplus B=\zz/M\zz$. Theorem \ref{divisors}
is proved. \qed

\bigskip
Our proof of Theorem \ref{three} will be based on the following
corollary of Theorem \ref{exponentials}.  

\begin{corollary}\label{constant}
Assume that $A\oplus B=\zz/M\zz$.  Define $A_m$ as in Theorem 
\ref{exponentials}. Let $N|M$, $c\in\zz\setminus B$, and
\[ b_m=b_m(c)=\#\{b\in B: \ (b-c,N)=m\}. \]
Then the quantity
\begin{equation}
 \sum_{m|N}\frac{b_m(c)A_m}{\phi(N/m)} \label{c.1}
\end{equation}
is independent of the choice of $c$.
\end{corollary}

{\bf Proof.}
Fix $N$ such that $N|M$. Apply Theorem \ref{exponentials} with $B$ replaced by
$C=B\cup\{c\}$, where $c\in\zz\setminus B$ will be allowed to
vary later on. Define $C(x)$ and ${\cal C}_d$ in the obvious way.

We first evaluate the terms on the right of (\ref{main.1}).  The term with $d=1$ is
$(|B|+1)^2 |A|^2/N$. If $d\neq 1,d|N$, then $d|M$, hence $\Phi_d(x)$ 
divides at least one of $A(x)$ or $B(x)$. If it divides $A(x)$, then $\cala_d=0$.  
If it divides $B(x)$, then for all roots $\xi$ of $\Phi_d(x)$ we have
$|C(\xi)|^2=|B(\xi)+\xi^c|^2=|\xi^c|^2=1$, hence ${\cal C}_d
=\#\{\xi:\Phi_d(\xi)=0\}=\phi(d)$.  Combining all this we find
that the right-hand side of (\ref{main.1}) equals
\begin{equation}
 \frac{1}{N}(|B|+1)^2|A|^2+\frac{1}{N}\sum_{d|N,d\neq 1}\cala_d. \label{main.right}
\end{equation}
Observe that this is independent of the choice of $c$.

Next, we have $C_m=B_m+2b_m$, hence the left side of (\ref{main.1}) equals
\begin{equation}
 \sum_{m|N}\frac{A_mB_m}{\phi(N/m)}+
\sum_{m|N}\frac{2b_mA_m}{\phi(N/m)} . \label{main.left}
\end{equation}
Comparing (\ref{main.left}) and (\ref{main.right}) we obtain (\ref{c.1})

We remark that (\ref{c.1}) can be computed explicitly if $N=M$.  Namely, choose
$c$ so that $c=b+kN\in B+N\zz$.  If $m|N,m\neq N$, then $b_m\neq
0$ implies that there is an $b'\in B$ such that $(c-b',M)=m$,
hence $(b-b',M)=m$ and $m\in D_B$.  By Theorem \ref{divisors},
$m\notin D_A$ and $A_m=0$. It follows that $b_mA_m=0$ for all
$m\neq M$. Moreover $b_M=1$ and $A_M=|A|$. Hence 
\begin{equation}
\sum_{m|N}\frac{b_m(c)A_m}{\phi(N/m)}=|A|\hbox{ if }N=M.\ \qed
\label{c.2} \end{equation}

\section{The tiling result}
\label{sec3}
\init


In this section we prove Theorem \ref{three}. Let $A,B,M$ be as in the statement
of the theorem.

Throughout the proof we will assume that $B$ is not contained in
$d\zz$ for any $d|M, d\neq 1$, for otherwise we may decompose
the tiling as in Lemma 2.5 of \cite{CM} and proceed by induction.
More precisely, suppose that the theorem is true for all sets $A'$ 
whose cardinality $|A'|$ divides, but is not equal to, $|A|$.  
Suppose further that $B\subset p\zz$.  From Lemma 2.5 of \cite{CM}
we have the decomposition
\[
A(x)=\sum_{i=0}^{p-1} x^{a_i}\bar A^{(i)}(x^p),
\]
where $A^{(i)}=\{a\in A:\ a\equiv i\ (\mod p)\}$, $a_i=\min(A^{(i)})$,
and $\bar A^{(i)}=\{a-a_i:\ a\in A^{(i)}\}/p$. Moreover, we have
\begin{equation}\label{x.e1}
|A^{(0)}|=|A^{(1)}|=\dots= |A^{(0)}|=|A|/p,
\end{equation}
\begin{equation}\label{x.e4}
\bar A^{(i)}\oplus p^{-1}B=\zz/Mp^{-1}\zz,
\end{equation}
\begin{equation}\label{x.e2}
S_{\bar A^{(0)}}=  S_{\bar A^{(1)}}=\dots=  S_{\bar A^{(p-1)}}
\end{equation}
and
\begin{equation}\label{x.e3}
S_A=\{p\}\cup S_{p\bar A^{(0)}}.
\end{equation}
Suppose that $\Phi_p(x),\Phi_q(x),\Phi_r(x)$ divide $A(x)$.  By 
(\ref{x.e3}), $\Phi_q(x)$ and $\Phi_r(x)$ divide $\bar A^{(0)}(x^p)$.
By Lemma 1.1(7) of \cite{CM}, $\Phi_q(x)$ and $\Phi_r(x)$ divide
$\bar A^{(0)}(x)$, hence also $\bar A^{(i)}(x)$ for all $i$ by (\ref{x.e2}).
Since $\bar A^{(i)}$ tiles $\zz$ by (\ref{x.e4}), it follows from the inductive
assumption that $\Phi_{qr}(x)$ divides $\bar A^{(i)}(x)$ for each $i$.
Using Lemma 1.1(7) of \cite{CM} again, we deduce that $\Phi_{pq}(x)$,
$\Phi_{pr}(x)$, $\Phi_{qr}(x)$, $\Phi_{pqr}(x)$ divide $\bar A^{(i)}(x^p)$ for each $i$,
hence they also divide $A(x)$.  Thus if we assume $A$ to be a set of
the smallest cardinality for which the theorem fails, the corresponding $B$
cannot be a subset of $p\zz$, $q\zz$, or $r\zz$.

The following notation will be used throughout this section.  We write
$[i,j,k]=n_{i,j,k}+pqr\zz$, where $n_{i,j,k}$ is the unique integer in
$\{0,1,\dots,pqr-1\}$ equal to $i (\mod p),j(\mod q),k(\mod
r)$.  We also write $[*,j,k]
=\bigcup_i [i,j,k]$, $[*,*,k]=\bigcup_{i,j}[i,j,k]$, etc. One can
think of the residues modulo $p,q,r$ as three-dimensional 
``coordinates", so that for example $[i,j,k]$ is a point, $[i,j,*]$ is a
vertical line, and $[*,*,k]$ is a horizontal plane.

\begin{lemma}
Let $B\subset\zz$. Assume that $0\in B$ and
\begin{equation}
B-B\subset p\zz\cup q\zz\cup r\zz. \label{pqr}
\end{equation}
Then at least one of the following holds:
\begin{equation}
 B\subset[*,j,k]\cup[i,*,k]\cup[i,j,*]\hbox{ for some }i,j,k,
\label{caseA}
\end{equation}
\begin{equation}
B\subset [0,0,0]\cup [i,j,0]\cup [i,0,k] \cup [0,j,k]\hbox{ for
some }i,j,k, \label{caseB}
\end{equation}
\begin{equation}
B\subset [*,*,0]\hbox{, or }B\subset [*,0,*],\hbox{ or }B\subset
[0,*,*]. \label{caseC}
\end{equation}
\label{three.l1}
\end{lemma}

\proof Suppose that (\ref{caseA}) and (\ref{caseC}) fail.  Then in
particular there is an $b\in B$ which is not in the set on the
right of (\ref{caseA}) with $i=j=k=0$, say $b\in [i,j,0]$ for some
$i,j\neq 0$. From our assumptions we have $B-b\subset p\zz\cup
q\zz\cup r\zz$. Hence
\[
\begin{array}{rl}B&\subset\Big([*,*,0]\cup[0,*,*]\cup[*,0,*]\Big)
\cap\Big([*,*,0]\cup[i,*,*]\cup[*,j,*]\Big)\\[3mm]
&=[*,*,0]\cup[i,0,*]\cup[0,j,*].
\end{array}\]
\relax From the failure of (\ref{caseC}) we get that at least one
of the following holds:
\smallskip

(a) there is a $b'\in B$ such that $b'\in [i,0,*]\setminus [*,*,0]$,
\smallskip

(b) there is a $b''\in B$ such that $b''\in [0,j,*]\setminus [*,*,0]$.
\smallskip

Suppose that (a) holds, then $B-b'\subset p\zz\cup q\zz\cup
r\zz$, hence $B\cap [0,*,*]\subset [*,0,0]\cup [i,*,0]$. Thus we
have (\ref{caseA}) unless (a) and (b) both hold.  In the latter
case, $b'\in [i,0,k]$ and $b''\in [0,j,k']$ for some
$k,k'\neq 0$. We then see from (\ref{pqr}) that $k=k'$ and that
(\ref{caseB}) holds. \qed

\bigskip
By Theorem \ref{sands}, at least one of the sets $D_A$, $D_B$
does not contain 1.  We deduce that at least one of $A-A$, $B-B$ satisfies
(\ref{pqr}), hence at least one of $A$, $B$ obeys the conclusions of
Lemma \ref{three.l1}. We will now show that $A$ cannot obey these
conclusions.  Indeed, we are assuming that the elements of $A$ are
distributed uniformly $\mod p$, $\mod q$, and $\mod r$. Hence each
plane $[i,*,*]$ contains exactly $|A|/p$ elements of $A$, etc.
This immediately contradicts (\ref{caseB}) and (\ref{caseC}),
since in both of these cases there are planes which do not contain
any elements of $A$.  Suppose now that (\ref{caseA}) holds. Assume
that $p<r$ and let $i'\neq i$.  By uniformity $\mod r$ and $\mod
p$, the planes $[*,*,k]$ and $[i',*,*]$ contain exactly $|A|/r$
and $|A|/p$ elements of $A$. But by (\ref{caseA}), all the
elements of $A$ which belong to $[i',*,*]$ are in fact in
$[i',j,k]$, hence in $[*,*,k]$.  This implies $|A|/p\leq |A|/r$,
which contradicts the assumption that $p<r$.

Thus $B$ satisfies one of (\ref{caseA}), (\ref{caseB}) (recall
that we assume that (\ref{caseC}) fails).

We record a simple lemma.

\begin{lemma}
Let $A\subset\zz$.  Then for any $m$ we have
\[ |\{(a,a')\in A\times A:\ m|a-a'\}|\geq \frac{|A|^2}{m}, \]
with equality if and only if the elements of $A$ are 
equi-distributed $\mod m$. \label{three.lemma2}
\end{lemma}

Let $N=pqr$. For $m|N$, we write
$\aa_m=\displaystyle{\frac{A_m}{\phi(N/m)}}$.  It suffices to prove that
\begin{equation}
\aa_m=\aa_{m'} \hbox{ for all }m,m'|N.\label{equivalence}
\end{equation}
Indeed, (\ref{equivalence}) implies that
\[ |A|^2=\sum_{m|pqr}A_m=\sum_{m|pqr}\phi(\frac{pqr}{m})A_{pqr}
=pqr A_{pqr},\]
and the theorem follows by Lemma \ref{three.lemma2}.

It remains to deduce (\ref{equivalence}) from Corollary \ref{constant}.   

\bigskip
{\bf Case 1a.} Assume that $B$ satisfies (\ref{caseB}) and that $p,q,r>2$.
We may then choose $I,J,K$
such that $I\neq 0,i$, $J\neq 0,j$, $K\neq 0,k$, and the planes
$[I,*,*], [*,J,*], [*,*,K]$ contain no elements of $B$. We first
compare (\ref{c.1}) with $c\in [I,J,K]$ and $c'\in [0,J,K]$. We then have
\[ b_1(c)=|B|, \ b_m(c)=0\hbox{ if }m\neq 1,\]
\[ b_1(c')=|B\cap[i,*,*]|,\ b_p(c')=|B\cap [0,*,*]|\neq 0,\
b_m(c')=0\hbox{ if }m\neq 1,p.\]
Substituting this in (\ref{c.1}) we see that
\[ |B|\aa_1=|B\cap[i,*,*]|\aa_1+|B\cap [0,*,*]|\aa_p,\]
hence $\aa_1=\aa_p$. Repeating this argument with $p$ replaced by $q$ and $r$, 
we obtain
\begin{equation}
\aa_p=\aa_q=\aa_r=\aa_1.  \label{sim2}
\end{equation}

With $I,J,K$ as above, let $c''\in[0,0,K]$, then
\[ b_{pq}=|B\cap[0,0,0]|\neq 0, \ b_{pr}=b_{qr}=b_{pqr}=0. \]
Comparing (\ref{c.1}) for $c$ and $c''$, and using also (\ref{sim2}), we 
find that $\aa_{pq}=\aa_1$.  Similarly for $\aa_{qr}$ and $\aa_{pr}$, hence
\begin{equation}
\aa_1=\aa_{pq}=\aa_{pr}=\aa_{qr}. \label{sim3}
\end{equation}

It only remains to prove that $\aa_{pqr}=\aa_1$.  But this follows by
applying (\ref{c.1}) and (\ref{sim2}), (\ref{sim3}) to $c$ as
above and $c'''\in[0,0,0]$.

\bigskip
{\bf Case 1b.} Assume now that $B$ satisfies (\ref{caseB}) and that $p=2$.
Let
\[ t=|B\cap[0,0,0]|,\ x=|B\cap[i,j,0]|,\ y=|B\cap[0,j,k]|,\ z=|B\cap[i,0,k]|.\]
Since the case when (\ref{caseA}) holds will be considered below, we may 
now assume that (\ref{caseA}) fails, and in particular that $t,x,y,z$ are all nonzero.
Choose $J,K$ such that $J\neq 0,j$, $K\neq 0,k$, and the planes
$[*,J,*], [*,*,K]$ contain no elements of $B$. We first
evaluate (\ref{c.1}) with $c\in [0,J,K],\ [0,J,0], [0,J,k]$, and find that
the following are all equal:
\[
\begin{array}{ll}
(t+y)\aa_2+(x+z)\aa_1&=C, \\[3mm]
t\aa_{2r}+y\aa_2+x\aa_r+z\aa_1&=C, \\[3mm]
y\aa_{2r}+t\aa_2+x\aa_r+z\aa_1&=C.
\end{array}\]
Therefore
\begin{equation}\label{b.e1}
\begin{array}{ll}
  t(\aa_2-\aa_{2r})+x(\aa_1-\aa_r)&=0,\\[3mm]
  y(\aa_2-\aa_{2r})+z(\aa_1-\aa_r)&=0.
  \end{array}
\end{equation}
Similarly, by considering (\ref{c.1}) with $c$ in $[1,J,K],\ [1,J,0],\ 
[1,J,k]$ we obtain that
\begin{equation}\label{b.e2}
\begin{array}{ll}
  x(\aa_2-\aa_{2r})+t(\aa_1-\aa_r)&=0,\\[3mm]
  z(\aa_2-\aa_{2r})+y(\aa_1-\aa_r)&=0.
  \end{array}
\end{equation}
Combining the first equations in (\ref{b.e1}), (\ref{b.e2}) we deduce that
$(x-t)(\aa_2-\aa_{2r}-\aa_1+\aa_r)=0$. Similarly, combining the second
equations we deduce that $(y-z)(\aa_2-\aa_{2r}-\aa_1+\aa_r)=0$.  It follows
that
\begin{equation}\label{b.e3}
  \aa_2-\aa_{2r}=\aa_1-\aa_r.
\end{equation}
Indeed, if (\ref{b.e3}) fails, we must have $x=t$ and $y=z$, in which case
$B$ is equi-distributed $\mod 2$ and $\Phi_2(\xi)$ divides both $A(\xi)$ and
$B(\xi)$.  This is easily seen to be impossible, e.g. by (T1).
We now substitute (\ref{b.e3}) in the first equation in (\ref{b.e2}):
\[(t+x)(\aa_2-\aa_{2r})=(t+x)(\aa_1-\aa_r)=0.\]
Since $t+x>0$, it follows that $\aa_1=\aa_r$ and $\aa_2=\aa_{2r}$.
We now repeat the same argument with $r$ replaced by $q$, and conclude that
\begin{equation}\label{b.e4}
  \aa_1=\aa_r=\aa_q,\ \aa_2=\aa_{2r}=\aa_{2q}.
\end{equation}

Next, we evaluate (\ref{c.1}) for $c$ in $[0,0,0],\ [1,0,0],\ [0,j,0],
\ [1,j,0],\ [0,0,k],\ [1,j,k]$.  Using also (\ref{b.e4}), we obtain that
\begin{equation}\label{b.e5}
  \begin{array}{lll}
  t\aa_{2qr}+x\aa_r+y\aa_2+z\aa_q&=t\aa_{2qr}+(x+z)\aa_1+y\aa_2&=C,\\[3mm]
  t\aa_{qr}+z\aa_{2q}+x\aa_{2r}+y\aa_1&=t\aa_{qr}+(x+z)\aa_2+y\aa_1&=C,\\[3mm]
  t\aa_{2r}+x\aa_{qr}+y\aa_{2q}+z\aa_1&=x\aa_{qr}+(y+t)\aa_2+z\aa_1&=C,\\[3mm]
  x\aa_{2qr}+t\aa_r+y\aa_q+z\aa_2&=x\aa_{2qr}+(y+t)\aa_1+z\aa_2&=C,\\[3mm]
  t\aa_{2q}+z\aa_{qr}+y\aa_{2r}+x\aa_1&=z\aa_{qr}+(t+y)\aa_2+x\aa_1&=C,\\[3mm]
  y\aa_{qr}+x\aa_{2q}+z\aa_{2r}+t\aa_1&=y\aa_{qr}+(x+z)\aa_2+t\aa_1&=C.  
  \end{array}
\end{equation}
\relax From equations 2,6 we have $(t-y)(\aa_{qr}-\aa_1)=0$, and from equations
3,5 $(x-z)(\aa_{qr}-\aa_1)=0$. Suppose first that $t\neq y$ or $x\neq z$, hence
$\aa_{qr}=\aa_1$.  Then we deduce from equations 2,4 that $\aa_2=\aa_{2qr}$. 
Substituting this in equations 1 and 2, we find that
\[ (x+z)\aa_1+(t+y)\aa_2=(x+z)\aa_2+(t+y)\aa_1,\]
hence $(x+z-t-y)(\aa_1-\aa_2)=0$.  Now $x+z\neq t+y$, since otherwise 
$B$ would be equi-distributed $\mod 2$ and we have already noted that this is
impossible. Therefore $\aa_1=\aa_2$, hence all the $\aa_m$ are equal and we are done. 

It remains to consider the case when $t=y$ and $x=z$.  Then we rewrite equations
1,2,3,5 in (\ref{b.e5}) as
\begin{equation}\label{b.e6}
  \begin{array}{lllll}
  2x\aa_1&+t\aa_2&&+t\aa_{2qr}&=C,\\[3mm]
  t\aa_1&+2x\aa_2&+t\aa_{qr}&&=C,\\[3mm]
  2t\aa_1&+x\aa_2&&+x\aa_{2qr}&=C,\\[3mm]
  x\aa_1&+2t\aa_2&+x\aa_{qr}&&=C.\\[3mm]
  \end{array}
\end{equation}
The determinant of the coefficient matrix is $-4(t^2-x^2)^2$.  If it were 0, we
would have $x=t=y=z$, and in particular $|B|=x+y+z+t=4t$ would be divisible by 4,
which contradicts the assumption that $|B|=2qr$.  Hence (\ref{b.e6}) has only
the trivial solution $\aa_1=\aa_2=\aa_{qr}=\aa_{2qr}$.  This together with
(\ref{b.e4}) implies that all the $\aa_m$ are equal, which completes the
proof for Case 1b.

\bigskip
{\bf Case 2.} Assume that $B$ satisfies (\ref{caseA}). Translating $B$ if 
necessary, we may assume that (\ref{caseA}) holds with $i=j=k=0$. Denote
\[ t=|B\cap [0,0,0]|, \]
\[ x_i=|B\cap [i,0,0]|,\ y_j=|B\cap [0,j,0]|,\
z_k=|B\cap [0,0,k]|,\ i,j,k>0,\]
\[ X=\sum x_i,\ Y=\sum y_i,\ Z=\sum z_i.\]
Since we are assuming that (\ref{caseC}) fails, we have $X,Y,Z\neq 0$.

Applying Corollary \ref{constant} to $c$ in $[i,j,k]$, $[i,j,0]$, $[i,0,k]$, $[0,j,k]$, 
$[i,0,0]$, $[0,j,0]$, $[0,0,k]$, $[0,0,0]$, where $i,j,k\neq 0$, we obtain 
that the following are all equal (denote the right-hand side by $C$):

\medskip
\begin{equation}\label{3e1}
\begin{array}{ll}
(X-x_i+Y-y_j+Z-z_k+t)\aa_1+x_i\aa_p+y_j\aa_q+z_k\aa_r&=C,
\\[3mm]
(X-x_i+Y-y_j+t)\aa_r+x_i\aa_{pr}+y_j\aa_{qr}+Z\aa_1&=C,
\\[3mm]
(X-x_i+Z-z_k+t)\aa_q+x_i\aa_{pq}+z_k\aa_{qr}+Y\aa_1&=C,
\\[3mm]
(Y-y_j+Z-z_k+t)\aa_p+y_j\aa_{pq}+z_k\aa_{pr}+X\aa_1&=C,
\\[3mm]
x_i\aa_{pqr}+(X-x_i+t)\aa_{qr}+(Y+Z)\aa_1&=C,
\\[3mm]
y_j\aa_{pqr}+(Y-y_j+t)\aa_{pr}+(X+Z)\aa_1&=C,
\\[3mm]
z_k\aa_{pqr}+(Z-z_k+t)\aa_{pq}+(X+Y)\aa_1&=C,
\\[3mm]
t\aa_{pqr}+X\aa_{qr}+Y\aa_{pr}+Z\aa_{pq}&=C.
\end{array}
\end{equation}
\medskip

We have to prove that this is possible if and only if all the
$\aa_m$ are equal.  We begin with a few lemmas.

\begin{lemma}
Let $A\subset\zz$, $|A|^2=pqrL$, $N=prq$.  Define $\aa_m$ as above. 
Assume that $\Phi_p,\Phi_q,\Phi_r,\Phi_{pr}$ divide $A(x)$. Then:
\begin{equation}\label{unif.e3}
  (q-1)\aa_{pr}+\aa_{pqr}=qL,
\end{equation}
\begin{equation}\label{unif.e4}
  (q-1)\aa_{r}+\aa_{qr}= (q-1)\aa_p+\aa_{pq}=qL,
\end{equation}
\begin{equation}\label{unif.e5}
  (q-1)\aa_{1}+\aa_{q}=qL.
\end{equation}
\label{three.l3}\end{lemma}

\proof 
We will first prove that if $\Phi_p,\Phi_q,\Phi_r$ divide $A(x)$, then:
\begin{equation}
\begin{array}{l}
(q-1)(r-1)\aa_p+(r-1)\aa_{pq}+(q-1)\aa_{pr}+\aa_{pqr}=qrL,
\\[3mm]
(p-1)(r-1)\aa_q+(r-1)\aa_{pq}+(p-1)\aa_{qr}+\aa_{pqr}=prL,
\\[3mm]
(p-1)(q-1)\aa_r+(q-1)\aa_{pr}+(p-1)\aa_{qr}+\aa_{pqr}=pqL,
\end{array}
\label{unif.e1}\end{equation}
and
\begin{equation}
\begin{array}{l}
(q-1)(r-1)\aa_1+(r-1)\aa_{q}+(q-1)\aa_{r}+\aa_{qr}=qrL,
\\[3mm]
(p-1)(r-1)\aa_1+(r-1)\aa_{p}+(p-1)\aa_{r}+\aa_{pr}=prL,
\\[3mm]
(p-1)(q-1)\aa_1+(q-1)\aa_{p}+(p-1)\aa_{q}+\aa_{pq}=pqL.
\end{array}
\label{unif.e2}\end{equation}
Indeed, from Lemma \ref{three.lemma2} with $m=p$ we have
\[ A_p+A_{pq}+A_{pr}+A_{pqr}=\frac{|A|^2}{p}=qrL,\]
and the first equation in (\ref{unif.e1}) follows by converting
the $A_m$ to $\aa_m$. Also, since $\sum_{m|pqr}A_m =|A|^2$,
from the displayed equation above we have
\[ A_1+A_q+A_r+A_{qr}=(1-\frac{1}{p})|A|^2=(p-1)qrL,\]
and the first equation in (\ref{unif.e2}) follows.  The remaining
equations in (\ref{unif.e1}), (\ref{unif.e2}) are similar.

Assume now that also $\Phi_{pr}(x)|A(x)$.  Applying Lemma
\ref{three.lemma2} with $m=pr$, we obtain
\[A_{pr}+A_{pqr}=\frac{|A|}{pr}=qL,\]
which implies (\ref{unif.e3}).  (\ref{unif.e4}) follows by
combining (\ref{unif.e3}) with the first and third equations in
(\ref{unif.e1}), and (\ref{unif.e5}) by combining the first
equation in (\ref{unif.e4}) with the first equation in
(\ref{unif.e2}). \qed

\begin{lemma} Suppose that $\aa_m$ solve (\ref{3e1}) with
$X,Y,Z\neq 0$, and that
\begin{equation} \aa_1=\aa_p,\ \aa_r=\aa_{pr},\ \aa_q=\aa_{pq},\
\aa_{qr}=\aa_{pqr}.\label{3e2}
\end{equation}
Then the $\aa_m$ are all equal. 
\label{three.l4}\end{lemma}

\proof Fix $i,j,k$. Plugging (\ref{3e2}) into (\ref{3e1}), we obtain
\medskip

\begin{equation}\label{3e3}
\begin{array}{ll}
(X+Y-y_j+Z-z_k+t)\aa_1+y_j\aa_q+z_k\aa_r&=C,
\\[3mm]
(X+Y-y_j+t)\aa_r+y_j\aa_{qr}+Z\aa_1&=C,
\\[3mm]
(X+Z-z_k+t)\aa_q+z_k\aa_{qr}+Y\aa_1&=C,
\\[3mm]
(Y-y_j+Z-z_k+t)\aa_p+y_j\aa_{q}+z_k\aa_{r}+X\aa_1&=C,
\\[3mm]
(X+t)\aa_{qr}+(Y+Z)\aa_1&=C,
\\[3mm]
y_j\aa_{qr}+(Y-y_j+t)\aa_{r}+(X+Z)\aa_1&=C,
\\[3mm]
z_k\aa_{qr}+(Z-z_k+t)\aa_{q}+(X+Y)\aa_1&=C,
\\[3mm]
(t+X)\aa_{qr}+Y\aa_{r}+Z\aa_{q}&=C.
\end{array}
\end{equation}
\medskip

From equations 2 and 6 in (\ref{3e3}) we have $X\aa_r=X\aa_1$, hence 
$\aa_r=\aa_1$.  Similarly, from equations 3 and 7 we have $X\aa_q=X\aa_1$, 
hence $\aa_q=\aa_1$.  We now have $\aa_1=\aa_p=\aa_q=\aa_r=\aa_{pr}=\aa_{pq}$.
Plugging this into equation 1 we obtain $(X+Y+Z+t)\aa_1=C$; this together with
equation 5 yields that $(X+t)\aa_{qr}=(X+t)\aa_1$, hence $\aa_{qr}=\aa_1$.  
By the last part of (\ref{3e2}) we also have $\aa_{pqr}=\aa_1$, which ends 
the proof. \qed

\bigskip
We now begin the proof of Theorem \ref{three} under the assumption that 
$B$ satisfies (\ref{caseA}). It suffices to consider the case when
\begin{equation}
 x_i=x,\ y_j=y,\ z_k=z \label{3e4}
\end{equation}
for some $x,y,z\neq 0$ and all $i,j,k\neq 0$.  (Hence $X=(p-1)x,
Y=(q-1)y, Z=(r-1)z$.)  Indeed, suppose for instance that $x_i\neq
x_{i'}$ for some $i,i'$.  Fix some $j,k$, and apply (\ref{3e1})
with $i,j,k$ and $i',j,k$.  From equations 1, 2, 3, 5 in
(\ref{3e1}) we find that (\ref{3e2}) holds, hence by Lemma
\ref{three.l4} all the $\aa_m$ are equal and we are done.

\begin{lemma}
Assume that $B$ satisfies (\ref{caseA}) and that (\ref{3e4})
holds.  Then:

\begin{itemize}
  \item $\Phi_{pq}(\xi)|B(\xi)$ if and only if $t=x+y+z-zr$;
  \item $\Phi_{qr}(\xi)|B(\xi)$ if and only if $t=x+y+z-xp$;
  \item $\Phi_{pr}(\xi)|B(\xi)$ if and only if $t=x+y+z-yq$;
  \item $\Phi_{pqr}(\xi)|B(\xi)$ if and only if $t=x+y+z$.
\end{itemize}
\label{three.l5}\end{lemma}

\proof  We have
\[ B(\xi)=t+x(\xi^{qr}+\xi^{2qr}+\dots+\xi^{(p-1)qr})
+y(\xi^{pr}+\xi^{2pr}+\dots+\xi^{(q-1)pr})\]
\[ +z(\xi^{pq}+\xi^{2pq}+\dots+\xi^{(r-1)pq})\]
\[=t+x(\Phi_p(\xi^{qr})-1)+y(\Phi_q(\xi^{pr})-1)
+z(\Phi_r(\xi^{pq})-1).\]
Hence
\[ B(e^{2\pi i/pq})=t+x\Phi_p(e^{2\pi ir/p})
+y\Phi_q(e^{2\pi ir/q}) +z\Phi_r(1)-x-y-z = t+zr-x-y-z,\]
and similarly
\[ \begin{array}{l}
 B(e^{2\pi i/qr})=t-x-y-z+px,\\[3mm]
B(e^{2\pi i/pr})=t-x-y-z+qy,\\[3mm]
B(e^{2\pi i/pqr})=t-x-y-z.
\end{array}\]
The lemma follows. \qed

\begin{corollary}
Let $B$ be as in Lemma \ref{three.l5}.
\begin{itemize}
  \item If $\Phi_{pqr}(\xi)|B(\xi)$, then none of
  $\Phi_{pq}(\xi),\Phi_{qr}(\xi),\Phi_{pr}(\xi)$ can divide $B(\xi)$.
  \item Assume that $|B|=pqr$, then at most one of
  $\Phi_{pq}(\xi),\Phi_{qr}(\xi),\Phi_{pr}(\xi)$ can divide $B(\xi)$.
\end{itemize}
\label{three.cor6}\end{corollary}

\proof The first part is obvious from Lemma \ref{three.l5}, since
$x,y,z\neq 0$.  Suppose now that $|B|=pqr$ and that
$\Phi_{pq},\Phi_{qr}$ divide $B(\xi)$.  By Lemma \ref{three.l5} we
have $t=x+y+z-zr=x+y+z-px$, hence $px=zr$, and in particular
$p|z$, $r|x$.  Moreover, adding up the elements of $B$ we obtain
\[|B|=pqr=t+(p-1)x+(q-1)y+(r-1)z=px+qy=qy+rz,\]
hence $qr|x$ and $pr|y$.  But then $\displaystyle{pqr=pqr
\frac{x}{qr}+pqr \frac{y}{pr}},$ therefore $x=0$ or $y=0$ -- a
contradiction.  \qed

\bigskip
We return to the proof of Theorem \ref{three}.  If
$\Phi_{pq}(x),\Phi_{qr}(x),\Phi_{pr}(x),\Phi_{pqr}(x)$ divide
$A(x)$, we are done.  Assume therefore that at least one of them
divides $B(x)$.  By Corollary \ref{three.cor6}, we only need to
consider two cases.

\bigskip
{\bf Case 2a:} $\Phi_{pq}(\xi)|B(\xi)$, $\Phi_{pr}(\xi)\Phi_{qr}(\xi)|A(\xi)$.
From Lemma \ref{three.l5} we have $t=x+y-Z$, which we substitute in (\ref{3e1}):

\medskip
\begin{equation}\label{3e6}
\begin{array}{ll}
(X+Y-z)\aa_1+x\aa_p+y\aa_q+z\aa_r&=C,
\\[3mm]
(X+Y-Z)\aa_r+x\aa_{pr}+y\aa_{qr}+Z\aa_1&=C,
\\[3mm]
(X+y-z)\aa_q+x\aa_{pq}+z\aa_{qr}+Y\aa_1&=C,
\\[3mm]
(Y+x-z)\aa_p+y\aa_{pq}+z\aa_{pr}+X\aa_1&=C,
\\[3mm]
x\aa_{pqr}+(X+y-Z)\aa_{qr}+(Y+Z)\aa_1&=C,
\\[3mm]
y\aa_{pqr}+(Y+x-Z)\aa_{pr}+(X+Z)\aa_1&=C,
\\[3mm]
z\aa_{pqr}+(x+y-z)\aa_{pq}+(X+Y)\aa_1&=C,
\\[3mm]
(x+y-Z)\aa_{pqr}+X\aa_{qr}+Y\aa_{pr}+Z\aa_{pq}&=C.
\end{array}
\end{equation}
\medskip

We also have from Lemma \ref{three.l3}:
\begin{equation}\label{3e7}
\begin{array}{l}
Y\aa_{pr}+y\aa_{pqr}=Y\aa_p+y\aa_{pq}=Y\aa_r+y\aa_{qr}
=Y\aa_1+y\aa_q=qyL,\\[3mm]
X\aa_{qr}+x\aa_{pqr}=X\aa_q+x\aa_{pq}=X\aa_r+x\aa_{pr}
=X\aa_1+x\aa_p=pxL,\\[3mm]
\end{array}
\end{equation}
where as before we denote $L=|A|^2/pqr$.  Plugging (\ref{3e7})
into (\ref{3e6}), we obtain:

\medskip
\begin{equation}\label{3e8}
\begin{array}{ll}
z(\aa_r-\aa_1)+pxL+qyL&=C,
\\[3mm]
Z(\aa_1-\aa_r)+pxL+qyL&=C,
\\[3mm]
z(\aa_{qr}-\aa_q)+pxL+qyL&=C,
\\[3mm]
z(\aa_{pr}-\aa_p)+pxL+qyL&=C,
\\[3mm]
(y-Z)\aa_{qr}+(Y+Z)\aa_1+pxL&=C,
\\[3mm]
(x-Z)\aa_{pr}+(X+Z)\aa_1+qyL&=C,
\\[3mm]
z\aa_{pqr}+(x+y-z)\aa_{pr}+(X+Y)\aa_1&=C,
\\[3mm]
Z(\aa_{pq}-\aa_{pqr})+pxL+qyL&=C.
\end{array}
\end{equation}
\medskip

\relax From equations 1,2 in (\ref{3e8}) we have $\aa_1=\aa_r$ and
$C=pxL+qyL$.  From equations 3,4,8 respectively we then have
$\aa_q=\aa_{qr}$, $\aa_p=\aa_{pr}$, $\aa_{pq}=\aa_{pqr}$.  Thus we
may apply Lemma \ref{three.l4} (with $p$ and $r$ interchanged) and
conclude that all the $\aa_m$ are equal.

\bigskip
{\bf Case 2b.} $\Phi_{pqr}(\xi)|B(\xi)$, $\Phi_{pq}(\xi)\Phi_{qr}(\xi)
\Phi_{pr}(\xi)|A(\xi)$.  By (\ref{unif.e3})--(\ref{unif.e5}) we have
\begin{equation}\label{3e11}
\begin{array}{l}
(p-1)\aa_{qr}+\aa_{pqr}=(p-1)\aa_q+\aa_{pq}=
(p-1)\aa_r+\aa_{pr}=(p-1)\aa_1+\aa_p=pL,\\[3mm]
(q-1)\aa_{pr}+\aa_{pqr}=(q-1)\aa_p+\aa_{pq}=
(q-1)\aa_r+\aa_{qr}=(q-1)\aa_1+\aa_q=qL,\\[3mm]
(r-1)\aa_{pq}+\aa_{pqr}=(r-1)\aa_r+\aa_{pr}=
(r-1)\aa_q+\aa_{qr}=(r-1)\aa_1+\aa_r=rL.
\end{array}\end{equation}
Thus we can compute all the $\aa_m$ if $\aa_1=\aa$ is given:
\begin{equation}\label{3e9}
\begin{array}{l}
\aa_p=pL-(p-1)\aa,\\[3mm]
\aa_q=qL-(q-1)\aa,\\[3mm]
\aa_r=rL-(r-1)\aa,\\[3mm]
\aa_{pq}=(p-1)(q-1)\aa-(pq-p-q)L,\\[3mm]
\aa_{pr}=(p-1)(r-1)\aa-(pr-p-r)L,\\[3mm]
\aa_{qr}=(q-1)(r-1)\aa-(qr-q-r)L,\\[3mm]
\aa_{pqr}=((p-1)(q-1)(r-1)+1)L-(p-1)(q-1)(r-1)\aa.
\end{array}\end{equation}
If $\Phi_{pqr}$ does not divide $A(\xi)$, by Lemma
\ref{three.lemma2} with $m=pqr$ we have $A_{pqr}=\aa_{pqr}>L$,
hence (from the last equation above) $L>\aa$.  We have to show
that this is impossible.

By Lemma \ref{three.l5} we have $t=x+y+z$.  We substitute this in the
last four equations in (\ref{3e1}):

\medskip
\begin{equation}\label{3e10}
\begin{array}{ll}
x\aa_{pqr}+(X+y+z)\aa_{qr}+(Y+Z)\aa_1&=C,
\\[3mm]
y\aa_{pqr}+(Y+x+z)\aa_{pr}+(X+Z)\aa_1&=C,
\\[3mm]
z\aa_{pqr}+(x+y+Z)\aa_{pq}+(X+Y)\aa_1&=C,
\\[3mm]
(x+y+z)\aa_{pqr}+X\aa_{qr}+Y\aa_{pr}+Z\aa_{pq}&=C.
\end{array}
\end{equation}
\medskip
(the remaining equations are equivalent).
We now plug in (\ref{3e11}).  From the last equation we have
\begin{equation}\label{3e13}
 xpL+yqL+zrL=C.
 \end{equation}
The remaining equations become
\begin{equation}\label{3e12}
\begin{array}{ll}
xpL+(y+z)\aa_{qr}+(Y+Z)\aa_1&=C,
\\[3mm]
yqL+(x+z)\aa_{pr}+(X+Z)\aa_1&=C,
\\[3mm]
zrL+(x+y)\aa_{pq}+(X+Y)\aa_1&=C.
\end{array}
\end{equation}
This adds up to
\[ xpL+yqL+zrL+(y+z)\aa_{qr}+(x+z)\aa_{pr}+(x+y)\aa_{pq}+2(X+Y+Z)\aa_1=3C,\]
hence by (\ref{3e13})
\[ (y+z)\aa_{qr}+(x+z)\aa_{pr}+(x+y)\aa_{pq}
= 2(pxL+qyL+rzL-X\aa_1-Y\aa_1-Z\aa_1).\]
By (\ref{3e11}), the left side equals 
\[ 2(pxL-X\aa_p+qyL-Y\aa_q+rzL-Z\aa_r).\]
But now we can use (\ref{3e9}). If $L>\aa$, we have
\[ \aa_p=pL-(p-1)\aa>\aa=\aa_1\]
and similarly $\aa_q>\aa_1,\aa_r>\aa_1$, which clearly contradicts
the above. \qed

\bigskip
{\bf Acknowledgement.} The first and third authors are supported in part by
the NSF. The second author is supported in part by NSERC.

{\sc Department of Mathematics, University of Georgia, Athens, GA 30602, U.S.A.}

{\it E-mail address:} {\tt andrew@math.uga.edu}
\medskip

{\sc Department of Mathematics, University of British Columbia, Vancouver, B.C. V6T 1Z2, Canada}

{\it E-mail address:} {\tt ilaba@math.ubc.ca}
\medskip

{\sc School of Mathematics, Georgia Institute of Technology, Atlanta, GA 30332, U.S.A.}

{\it E-mail address:} {\tt wang@math.gatech.edu}
\medskip

\end{document}